\documentclass[11pt]{article}
\usepackage[final]{epsfig}
\usepackage{graphics}
\usepackage{amsmath}
\usepackage{amsfonts}
\usepackage{latexsym}
\usepackage{amssymb}
\usepackage{graphicx}
\usepackage{url}
\usepackage{epstopdf}
\usepackage{hyperref}
\usepackage{color}

\newtheorem{lemma}{Lemma}

\newtheorem{theorem}{Theorem}

\newtheorem{corollary}[lemma]{Corollary}
\newtheorem{conjecture}[lemma]{Conjecture}

\newcommand{\g}{{\gamma}}
\newcommand{\G}{{\Gamma}}

\newcommand{\proofend}{$\Box$\bigskip}

\newcommand{\R}{{\mathbb R}}
\newcommand{\Z}{{\mathbb Z}}
\newcommand{\RP}{{\mathbb {RP}}}

\def\proof{\paragraph{Proof.}}

\begin{document}

\title{Iterating skew evolutes and skew involutes: a linear analog of the bicycle kinematics}

\author{
 Serge Tabachnikov\footnote{
Department of Mathematics,
Pennsylvania State University,
University Park, PA 16802,
USA;
tabachni@math.psu.edu}
} 

\date{\today}

\maketitle

\begin{abstract}
The evolute of a plane curve is the envelope of its normals. Replacing the normals by the lines that make a fixed angle with the curve yields 
a new curve, called the evolutoid. We prefer the term ``skew evolute", and we study the geometry and dynamics of the skew evolute map and of its inverse, the skew involute map. The relation between a curve and its skew evolute is analogous to the relation between the rear and front bicycle tracks, and this connections with the bicycle kinematics (a considerably more complicated subject) is our motivation for this study.
\end{abstract}

\paragraph{Introduction and motivation.}
The evolute of a smooth plane curve is the envelope of its normals. In this article we consider the following modification of this construction.

\begin{figure}[ht]
\centering
\includegraphics[width=.3\textwidth]{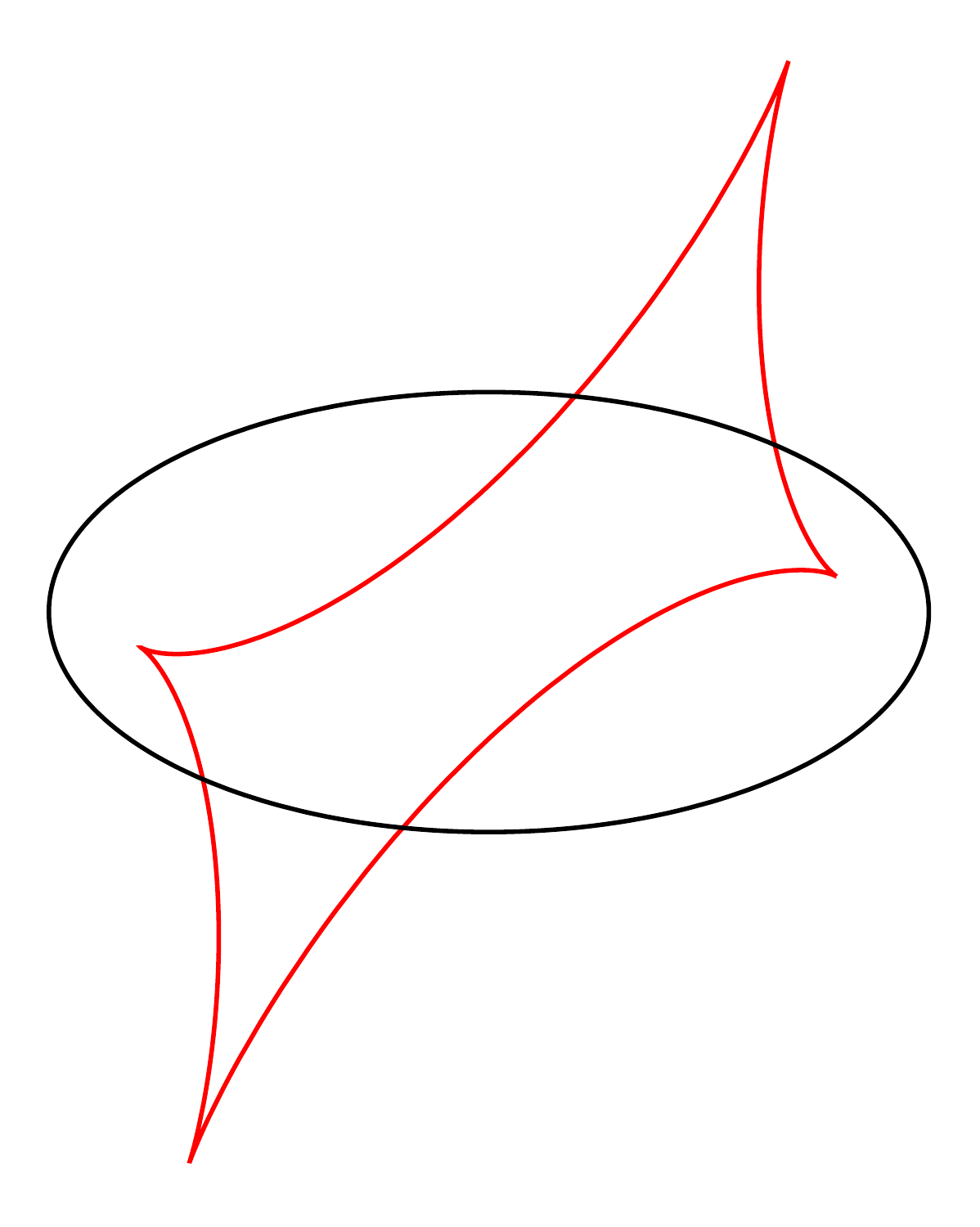}
\caption{A skew evolute of an ellipse.}	
\label{evellipse}
\end{figure}

Let $\g$ be a smooth oriented curve and $\alpha$ a fixed angle. Turn each tangent line of $\g$ through angle $\alpha$ about the tangency point, and let $\G$ be the envelope of this 1-parameter family of lines. We call $\G$ a {\it skew evolute} of $\g$ and write  $\G={\mathcal E}_\alpha(\g)$. See Figure \ref{evellipse}.
The usual evolute corresponds to the case $\alpha=\pi/2$, and if $\alpha=0$, then $\G=\g$. Likewise, we call $\g$ a {\it skew involute} of $\G$ and write $\g={\mathcal I}_\alpha (\G)$. 

This subject  goes back to the early 18th century, see \cite{Re} (we learned about this reference from \cite{GW}). However, it continues to attract attention; see \cite{AA,AM,GW,Ha,JC,RGS,Z} for a sampler of this century work.

What we called ``skew evolute" is traditionally called ``evolutoid". The reason we use a different term is to emphasize the similarity with the classical evolutes and involutes. Indeed, what we call ``skew involutes" were called ``tanvolutes" in \cite{AM}. It seems that the terminology has not  completely crystalized yet.

A study of skew evolutes necessarily involves curves with cusps; indeed, the evolute of a closed simple curve has at least four cusps, as the classical 4-vertex theorem implies. We study the dynamics of the transformations ${\mathcal E}_\alpha$ and ${\mathcal I}_\alpha$ on the class of curves called {\it hedgehogs}.

\begin{figure}[ht]
\centering
\includegraphics[width=.4\textwidth]{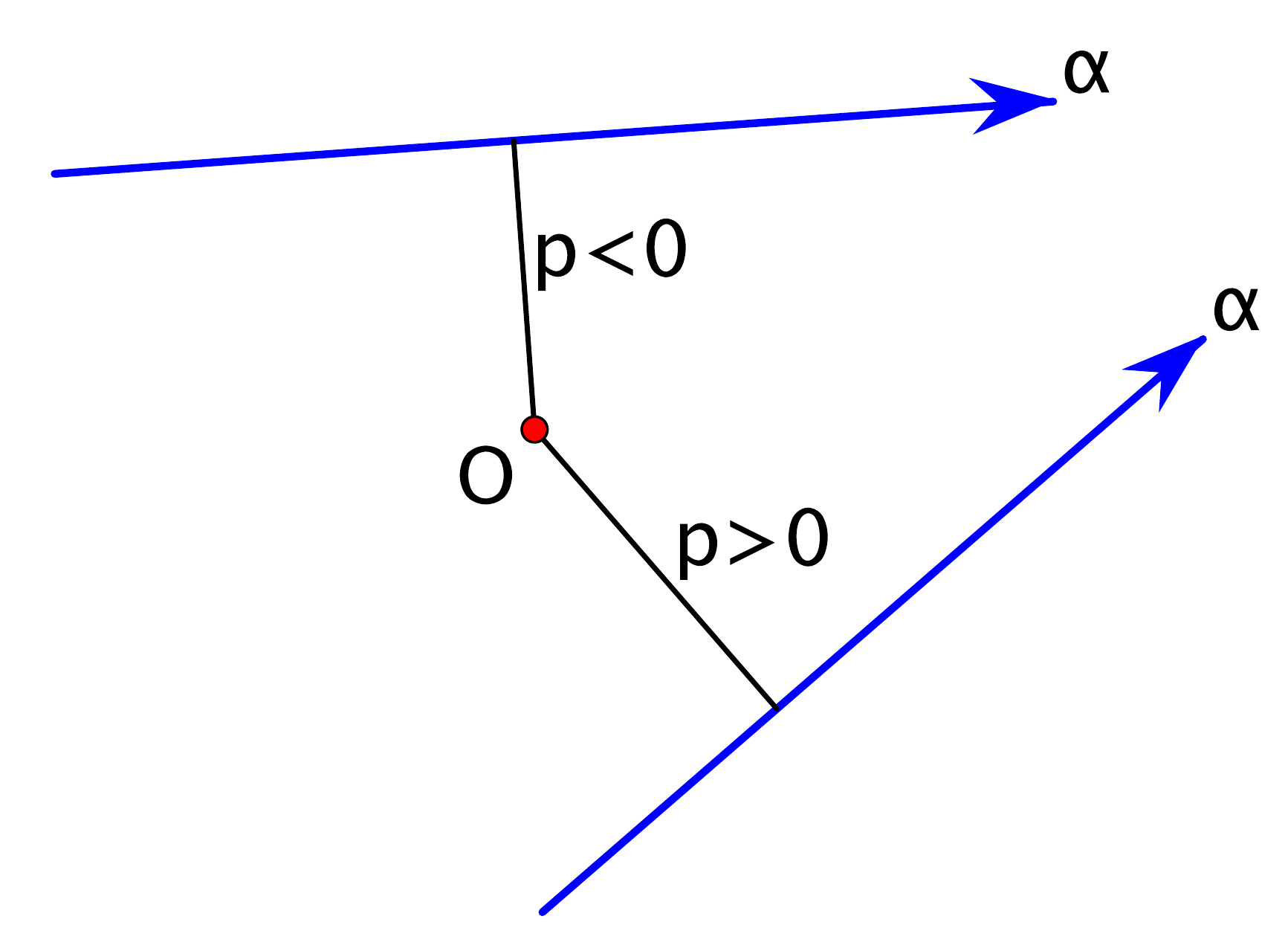}
\caption{The space of (co)oriented lines.}	
\label{lines}
\end{figure}

An oriented line in the plane is characterized by its direction $\alpha$ and the support number $p$, the signed distance from the origin to the line, Figure \ref{lines}. The coorientation of an oriented line is given by the direction $\phi=\alpha-\pi/2$.

Let $\g$ be an oriented smooth strictly convex closed curve. It can be parameterized by $\phi\in S^1=\R/2\pi\Z$, the direction of its outward normal vectors, and the support numbers of the tangent lines are given by a function $p(\phi)$. This is the support function of $\g$. 
The support function uniquely characterizes the curve and is determined by the curve, except that a change  of the origin amounts to adding to $p(\phi)$ a first harmonic, a linear combination of $\cos\phi$ and $\sin\phi$. 

In particular, the equation of the curve, defined by its support function, is
\begin{equation*} \label{eq:curv}
\g(\phi)=(p(\phi)\cos\phi-p'(\phi)\sin\phi,p(\phi)\sin\phi+p'(\phi)\cos\phi).
\end{equation*}
The length of $\g$ and the area bounded by it are given by
$$
L=\int_0^{2\pi} p(\phi)\ d\phi,\ \ A=\frac{1}{2} \int_0^{2\pi} [p^2(\phi)-(p')^2(\phi)]\ d\phi,
$$
and the curvature radius of $\g$ by $p(\phi)+p''(\phi)$, see, e.g., \cite{Sa}.

Replacing a curve by its equidistant curve amounts to adding a constant to the support function. The curves, that are equidistant to convex ones, still do not have inflections and are characterized by their support functions $p:S^1\to\R$, but they may have cusps, where the radius of curvature vanishes. The tangent lines are well defined at the cusps, and the coorientation is continuous therein (unlike the orientation that is reversed at the cusps).

The cooriented curves described by the support functions are  called {\it hedgehogs}. The orientation of a smooth arc of a hedgehog is obtained from its coorientation by a $90^\circ$ rotation in the positive direction.
An equivalent characterization of hedgehogs is that they are equidistant to convex curves. The above formulas for the length and area are still valid, but these quantities are signed (for example, the sign of the length changes as one traverses a cusp).

\begin{figure}[ht]
\centering
\includegraphics[width=.65\textwidth]{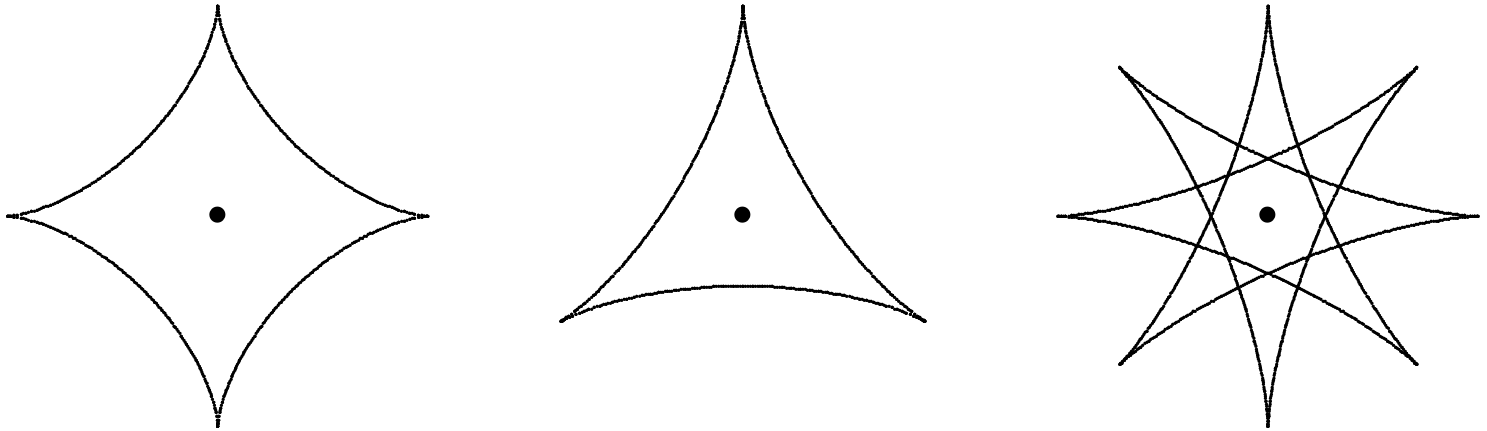}
\caption{Hypocycloids of order $2,3$, and $4$ (the middle curve is traversed twice).}	
\label{hypo}
\end{figure}

A {\it hypocycloid} is a hedgehog whose support function is a pure harmonic, a linear combination of $\cos(k\phi)$ and $\sin(k\phi)$. The number $k\ge 2$ is the order of a hypocycloid, see Figure \ref{hypo}. We consider circles as the hypocycloids of order zero.

There are two motivations for this study. First, this is a generalization of the work done in \cite{AFITT}, where iterations of evolutes and involutes were considered, both in the continuous and discrete settings (when curves are replaced by polygons).  

Second, there is a relation with the recent study of bicycle kinematics that we now describe. 

\begin{figure}[ht]
\centering
\includegraphics[width=1\textwidth]{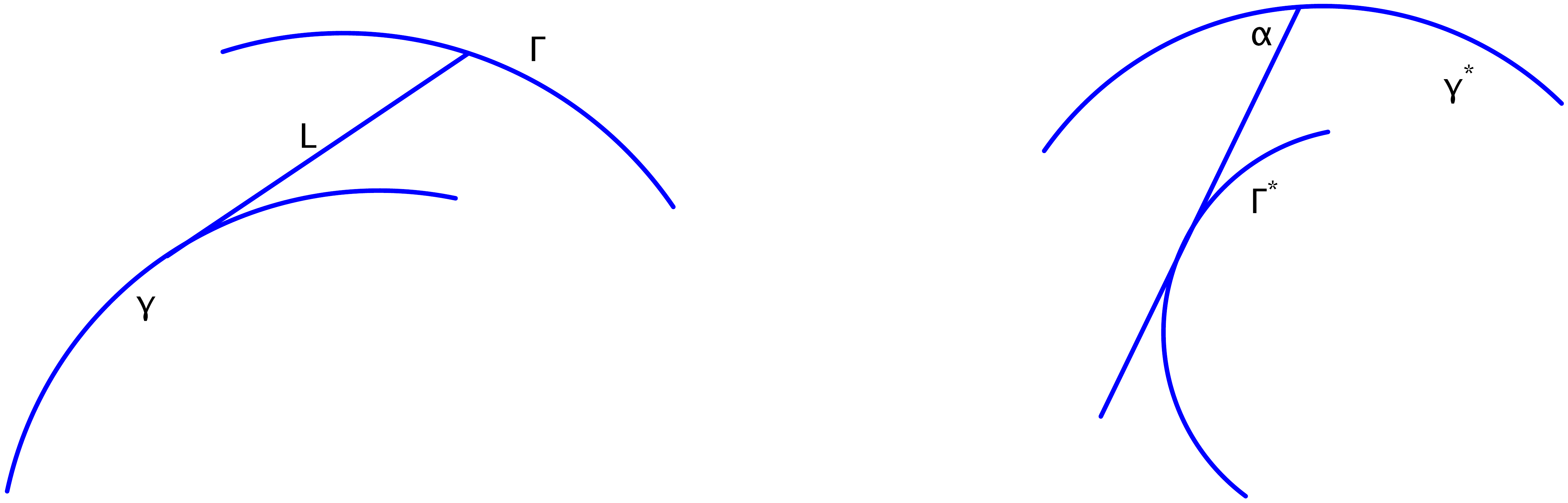}
\caption{The correspondence between the rear and front bicycle tracks and the dual  picture.}	
\label{db}
\end{figure}

Bicycle is modeled as an oriented segment of fixed length $L$ that can move in such a way that the velocity of its rear end is always aligned with the segment (the rear wheel is fixed on the bicycle frame). The bicycle leaves two tracks, the rear track $\g$ and the front track $\G$, and they are related as shown on the left of Figure \ref{db}. See, e.g., \cite{BLPT,FLT,LT}. 

This model of bicycle can be also considered in the spherical geometry, see \cite{HPZ}. In the spherical geometry one has a duality between points and oriented great circles, the pole-equator correspondence. This spherical duality extends to smooth curves and, applied to the left part of Figure \ref{db}, it yields the right part, where the angle $\alpha$ equals the spherical length $L$. 

However, we consider the right part of Figure \ref{db} as drawn in the plane. In this way, the map that takes the rear bicycle track to the front track is analogous to the map that takes a curve to its skew evolute. Unlike the former map, the latter one is  linear, and it is much easier to study.  In particular, as we shall see, the bicycle analogs of some results concerning skew evolutes and skew involutes remain open problems.

\paragraph{Known results.} In this section we present known results on skew evolutes and skew involutes.

Let $\g$ be a hedgehog with the support function $p(\phi)$, and let $\G={\mathcal E}_\alpha$ be its skew-evolute with the support function $q(\phi)$. Then
\begin{equation} \label{eq:qtop}
q(\phi)=p(\phi-\alpha)\cos\alpha+p'(\phi-\alpha)\sin\alpha,
\end{equation}
see \cite{JC}.
In particular, $L(\G)=\cos\alpha\ L(\g)$.

Denote the linear differential operator on the right hand side of (\ref{eq:qtop}) by ${\mathcal D_\alpha} (p)(\phi)$.

The {\it Steiner point}, or the {\it curvature centroid}, $St(\g)$, of a curve $\g$ is its center of mass with the density equal to the curvature. In terms of the support function, it is given by
$$
St(\g)=\frac{1}{\pi} \int_0^{2\pi} p(\phi)\ (\cos\phi,\sin\phi)\ d\phi.
$$

A hedgehog $\g$ and its skew evolute $\G={\mathcal E}_\alpha(\g)$ share their Steiner points, see \cite{AA}.

For a quick proof, note that the Steiner point is characterized by the condition that, if it is chosen as the origin, then the support function is $L^2$-orthogonal to the first harmonics, that is, its Fourier expansion does not contain the first harmonics. This property is preserved by the operator ${\mathcal D}_\alpha$, and the result follows.

The evolute of a curve is the locus of the centers of its osculating circles. For skew evolutes, the role of circles is played by the logarithmic spirals. 

A logarithmic spiral centered at the origin is characterized by the property that the position vector of every point makes a constant angle $\alpha$ with the direction of the curve at this point. If $\alpha=\pi/2$, the spiral is a circle.

Call such logarithmic spirals $\alpha$-spirals. They form a 1-parameter family of curves.  Allowing parallel translation of the origin, results in a 3-parameter family of $\alpha$-spiral (similarly to circles). It follows that, for every $\alpha$, a smooth curve has an osculating $\alpha$-spiral at every point (it  approximates the curve to second order). A hyper-osculating $\alpha$-spiral is tangent to the curve to higher order.

Therefore the skew evolute ${\mathcal E}_\alpha (\g)$ is the locus of centers of the osculating $\alpha$-spirals of the curve $\g$. The cusps of ${\mathcal E}_\alpha (\g)$ correspond to hyper-osculating $\alpha$-spirals, see \cite{Wu}.

The cusps of a skew evolute happen when its radius of curvature $r(\phi)$ vanishes. In view of equation (\ref{eq:qtop}), this amounts to the equation $r \cos\alpha + r' \sin\alpha=0$, or $r'/r=-\cot\alpha$. See \cite{GW} for a study of cusps of skew evolutes.

Let $\G$ be a hedgehog. Given $\alpha$, does $\G$ have a closed skew-involute, and if so, how many? For $\alpha=\pi/2$, the involute is provided by the string construction, and a necessary and sufficient condition for it to close up is that the signed length of $\G$ vanishes, in which case one has a 1-parameter family of involutes. 

However, if $\alpha \neq \pi/2$, then there exists a unique closed skew-involute ${\mathcal I}_\alpha (\G)$, see \cite{AM}. The reason is that the monodromy of the linear differential equation (\ref{eq:qtop}) is a homothety of the real line with coefficient $\ne 1$. Such a map has a unique fixed point, corresponding to the desired periodic solution.

Comparing with the bicycle kinematics, we observe the following difference. Given a closed front bicycle track, the rear track is  determined by a first order ordinary Riccati differential equation, see, e.g., \cite{FLT}. Unlike the case of skew involutes, the monodromy of this equation takes values in the group $SL(2,\R)$, acting on the circle $\RP^1$ of the initial positions of the bicycle by fractional-linear transformations. Generically, such a transformation has none  or two fixed points.

\paragraph{Three examples.} The following examples concern locally convex curves that are not closed, and their support functions are not periodic anymore. However formula (\ref{eq:qtop}) is still valid.
\medskip

\begin{figure}[ht]
\centering
\includegraphics[width=.55\textwidth]{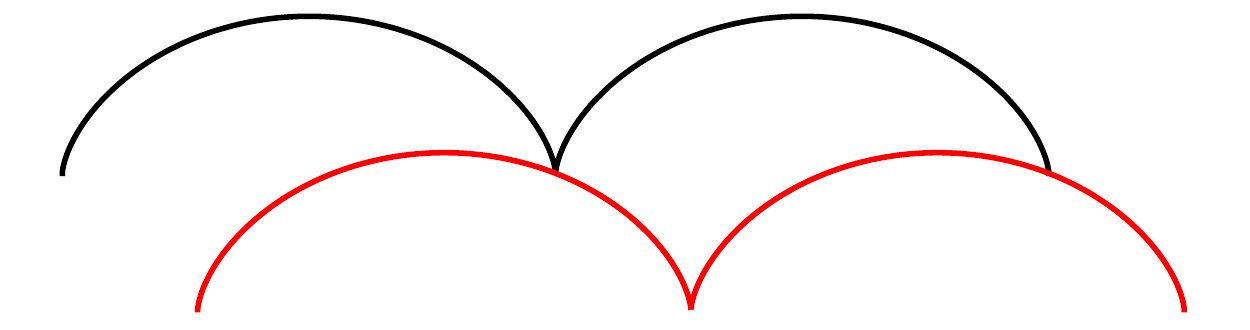}\quad
\includegraphics[width=.4\textwidth]{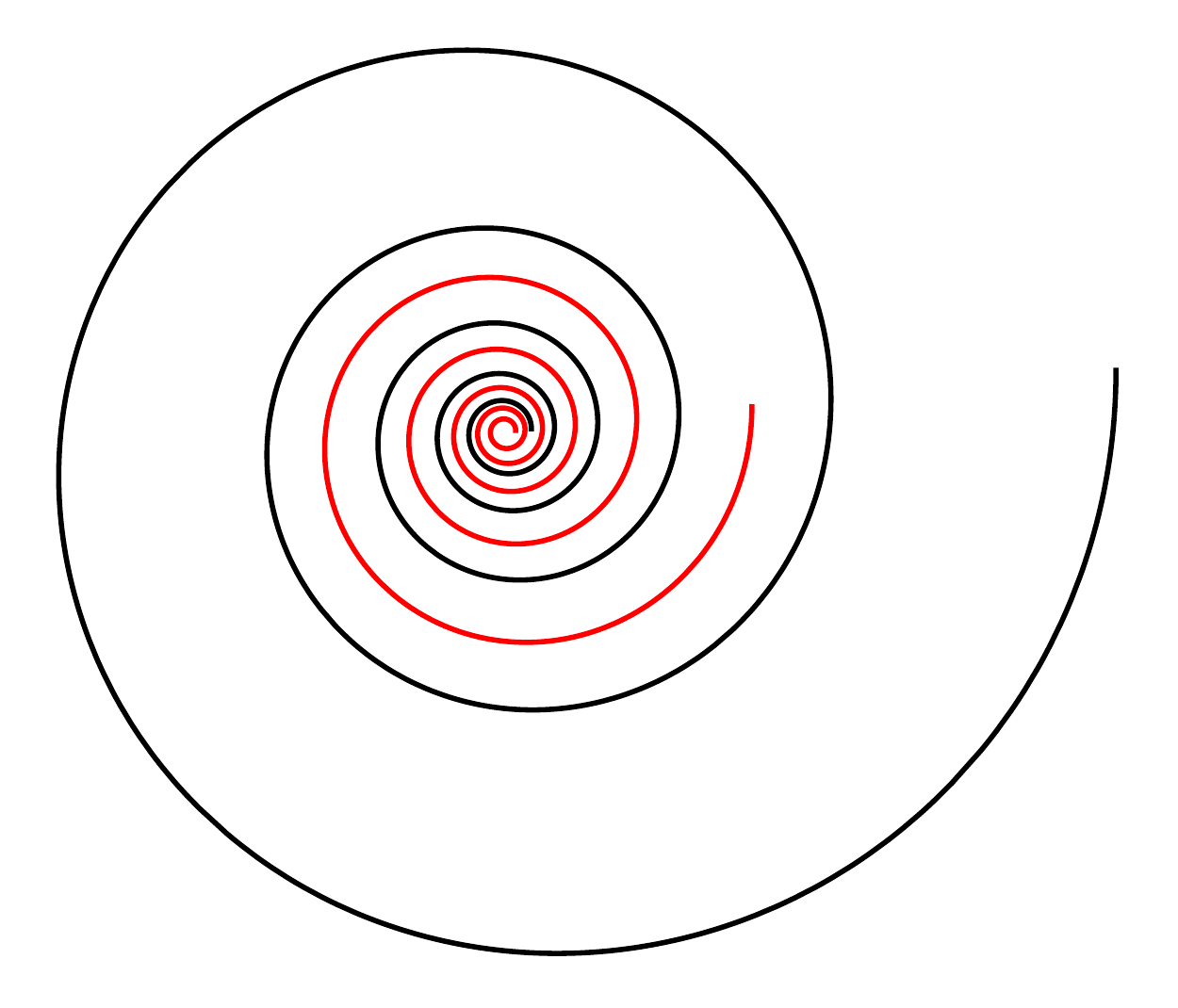}
\caption{Left: a cycloid and its skew evolute. Right: a logarithmic spiral and its skew evolute.}	
\label{cycloid}
\end{figure}

{\bf Cycloid.}
It is well known that the evolute of a cycloid is congruent to the cycloid by parallel translation. The same holds for skew evolutes, see Figure \ref{cycloid}. 

Indeed, the support function of a cycloid is $p(\phi)=-\phi \cos\phi$. Using equation (\ref{eq:qtop}), we find that the support function of the skew evolute is
$$
q(\phi)=-\phi \cos\phi + (\alpha-\cos\alpha \sin\alpha)\cos\phi + \sin^2\alpha \sin\phi.
$$
Thus the support function has changed by a first harmonic, which amounts to a parallel translation of the curve.

\medskip

{\bf Logarithmic spiral.}
Logarithmic spirals are  congruent to their skew evolutes by rotation, see Figure \ref{cycloid}. Indeed, the support function of a logarithmic spiral is $p(\phi)=e^{c\phi}$. Hence the support function of its skew evolute is
$$
q(\phi)=e^{-c\alpha} (\cos\alpha + c \sin\alpha) e^{c\phi},
$$
which is obtained from $e^{c\phi}$ by a parameter shift. If $c=-\cot\alpha$, the skew evolute reduces to a point.

A slight generalization is a curve $\g$ whose support function is $p(\phi)= c_1 e^{b_1 \phi}+c_2 e^{b_2 \phi}$. If 
$$
(\cos\alpha+b_1\sin\alpha)^{b_2} = (\cos\alpha+b_2\sin\alpha)^{b_1},
$$
then ${\mathcal E}_\alpha (\g)$ is congruent to $\g$ by  rotation. 
\medskip

\begin{figure}[ht]
\centering
\includegraphics[width=.6\textwidth]{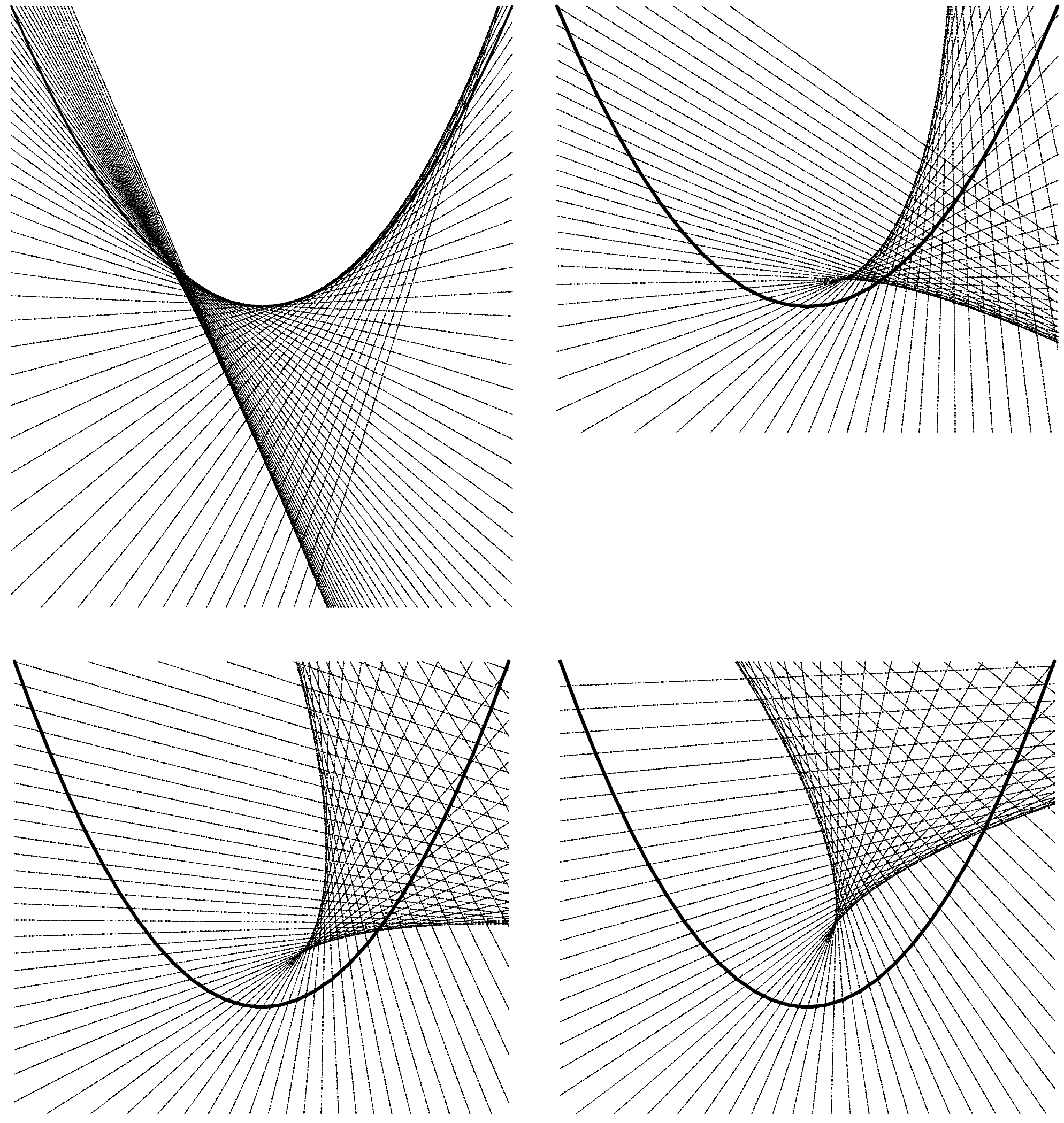}
\caption{Skew evolutes of a parabola with $\alpha=\pi/10, \pi/5, 3\pi/10$, and $2\pi/5$.}	
\label{parabola}
\end{figure}

{\bf Parabola.} A calculation, that we do not present, shows that the skew evolute of the parabola $(t,t^2/2)$ has a cusp for $3t=-\cot\alpha$, see Figure \ref{parabola}. Thus the skew evolute of a parabola  has a unique cusp for every $\alpha \in (0,\pi)$.

\paragraph{New results.} Now we present  results that, to the best of our knowledge, are not found  in the literature. 

First, let us look at the above defined linear operator ${\mathcal D_\alpha}$ in detail. It preserves the 2-dimensional space of $k$th harmonics. In the basis $(\cos(k\phi),\sin(k\phi))$, it is given by the matrix
\begin{equation} \label{eq:mat}
\begin{pmatrix}
\cos^2\alpha+k\sin^2\alpha,&(k-1)\cos\alpha\sin\alpha\\
-(k-1)\cos\alpha\sin\alpha,&\cos^2\alpha+k\sin^2\alpha
\end{pmatrix}.
\end{equation}
This is a linear similarity, a composition of rotation and dilation; the dilation coefficient is equal to $\sqrt{1+(k^2-1)\sin^2\alpha}$.
In particular, for $\alpha\ne \pi/2$, the operator ${\mathcal D_\alpha}$ is invertible.

Since similarities with a fixed center commute, we have

\begin{corollary} \label{cor:commute}
One has 
$$
{\mathcal E}_\alpha \circ {\mathcal E}_\beta  = {\mathcal E}_\beta \circ {\mathcal E}_\alpha\ \ {\rm and}\ \ 
{\mathcal I}_\alpha \circ {\mathcal I}_\beta  = {\mathcal I}_\beta \circ {\mathcal I}_\alpha.
$$
\end{corollary}

Next we ask how  the shape of a hedgehog evolves under iterations of the skew evolute or skew involute operations. 

\begin{theorem} \label{thm:iter}
(i) Assume that the support function of $\g$ is a  trigonometric polynomial of degree $d$. Then the iterated skew evolutes of $\g$ converge, in shape, to a hypocycloid of order $d$.\\
(ii) If the Fourier series of the support function of $\g$ has a free term, then its iterated skew involutes converge, in shape, to a circle. If the Fourier series starts with $d$th harmonics, then the  iterated skew involutes converge, in shape, to a hypocycloid of order $d$.
\end{theorem}

\proof
For the first statement, formula (\ref{eq:mat}) implies that, under iterations, the highest harmonics grow faster than the lower ones. This implies the result.

Likewise, under ${\mathcal D_\alpha}^{-1}$, the free term of the Fourier series is multiplied by $1/\cos\alpha > 1$, whereas the space of $k$th harmonics is stretched by the factor $1/\sqrt{1+(k^2-1)\sin^2\alpha}<1$. In the first case,  the free term dominates under iterations, and in the second case, so does the first non-trivial harmonic.
\proofend

\begin{corollary} \label{cor:sim}
A hedgehog is similar to its skew involute if and only if it is a hypocycloid.
\end{corollary}

In the case of evolutes and involutes ($\alpha=\pi/2$),  the above results were obtained in \cite{AFITT}. The next theorem extends another result  in \cite{AFITT} from evolutes to skew evolutes.

\begin{theorem} \label{thm:cusps}
Assume that the support function $p(\phi)$ of a hedgehog $\g$ is not a trigonometric polynomial, that is, its Fourier expansion contains infinitely many terms. Then the number of cusps of the iterated skew evolutes increases without bound.
\end{theorem}

\proof
The proof consists of two steps. 
\smallskip

Claim 1: {\it The number of sign changes of the functions ${\mathcal D}_\alpha^n(p)$ increases without bound as $n\to\infty$.} 

This is a slight generalization of the theorem by Polya and Wiener \cite{PW} where the case of the operator $p \mapsto p'$ is considered. 
Since this argument is not sufficiently well known, we present it here.

Let 
$$
p(\phi)=\sum_{k\in\Z} a_k e^{ki\phi},\ a_{-k}=\bar a_k,
$$
be the Fourier expansion of $p$. It suffices to proof the statement for a simpler operator ${\mathcal F}(p)=p'+cp$, that is,
$$
{\mathcal F}: p \mapsto \sum_{k\in\Z} (c+ik) a_k e^{ki\phi}.
$$
The claim is that if $a_m \neq 0$ then, for   sufficiently large $n$, the function ${\mathcal F}^n(p)$ has at least $2m$ sign changes.

Let ${\mathcal Z}(f)$ denote the number of sign changes of a periodic function $f$.
A version of Rolle's theorem, Lemma 1 in \cite{PW}, asserts that, for every $b\in\R$,
$$
{\mathcal Z}\left(\sum_{k\in\Z} a_k e^{ki\phi}\right) \ge {\mathcal Z}\left(\sum_{k\in\Z} \frac{a_k e^{ki\phi}}{b^2+k^2}\right).
$$
Apply this to $f={\mathcal F}^n(p)$, $b^2=m^2+2c^2$, and iterate the inequality $n$ times:
\begin{equation*}
\begin{split}
{\mathcal Z}\left(\sum_{k\in\Z} (c+ik)^n a_k e^{ki\phi}\right) \ge &{\mathcal Z}\left(\sum_{k\in\Z} \frac{(c+ik)^n a_k e^{ki\phi}}{(m^2+2c^2+k^2)^n}\right) = \\&{\mathcal Z}\left(\sum_{k\in\Z} \frac{2[\sqrt{c^2+m^2} (c+ik)]^n a_k e^{ki\phi}}{(m^2+2c^2+k^2)^n}\right).
\end{split}
\end{equation*}

Let $q_n(\phi)=\sum_{k\in\Z} c_k e^{ki\phi}$ be the function on the right. Then
$$
|c_k|= \left(\frac{2\sqrt{c^2+m^2}\sqrt{c^2+k^2}}{m^2+2c^2+k^2}  \right)^n |a_k|.
$$
One has
$$
\frac{2\sqrt{c^2+m^2}\sqrt{c^2+k^2}}{m^2+2c^2+k^2} < 1,
$$
unless $k=m$, in which case this coefficient equals 1. This implies that, for sufficiently large $n$, 
$$
|c_m| > \sum_{k\neq m} |c_k|,
$$
as in \cite{PW}. For such $n$,  ${\mathcal Z}(q_n)$ equals the number of sign changes of its $m$th harmonic, that is, equals $2m$, as needed.
\smallskip

Claim 2: {\it If the support function of a hedgehog $\g$ has $2m$ sign changes, then $\g$ has at least $m$ cusps.} 

Indeed, if the support function of $\g$ has $2m$ zeros, then there are $2m$ tangents from the origin  $O$ to $\g$. Each arc of a hedgehog between its cusps is convex, and there are at most two tangents from $O$ to it. Therefore there must be at least $m$ such arcs, and at least as many cusps. 
\proofend

Theorems \ref{thm:iter} and \ref{thm:cusps} imply 

\begin{corollary} \label{cor:cusps}
If all iterated skew evolutes of a hedgehog $\g$ are free from cusps, then $\g$ is a circle.
\end{corollary}

What is an analog of this statement in terms of the bicycle model? Since the projective duality interchanges cusps and inflections, we are led to the following formulation.

Let $\g$ be a smooth oriented closed curve, $L$ a fixed positive number. Denote by ${\mathcal T}(\g)$ the locus of endpoints of the positve tangent segments to $\g$ of length $L$. That is, ${\mathcal T}(\g)$ is the front track of the bicycle whose rear track is $\g$.

\begin{conjecture}
Assume that all iterations ${\mathcal T}^k(\g), k\ge 0,$ are convex curves. Then $\g$ is a circle. 
\end{conjecture} 

\paragraph{An ``integrable" map on hedgehogs.}
 Given a  bicycle rear track, one can traverse it in the opposite directions, creating two front tracks. This relation between curves is completely integrable, see \cite{BLPT}. Equivalently, two smooth curves, $\G_1$ and $\G_2$, are in the bicycle correspondence if two points, $x_1$ and $x_2$, can traverse them in such a way that the distance between them remains constant (twice the bicycle frame) and the velocity of the midpoint of the segment $x_1x_2$ is aligned with this segment (that is how the read wheel moves). 

An analog of this relation in our setting is as follows.

\begin{figure}[ht]
\centering
\includegraphics[width=.45\textwidth]{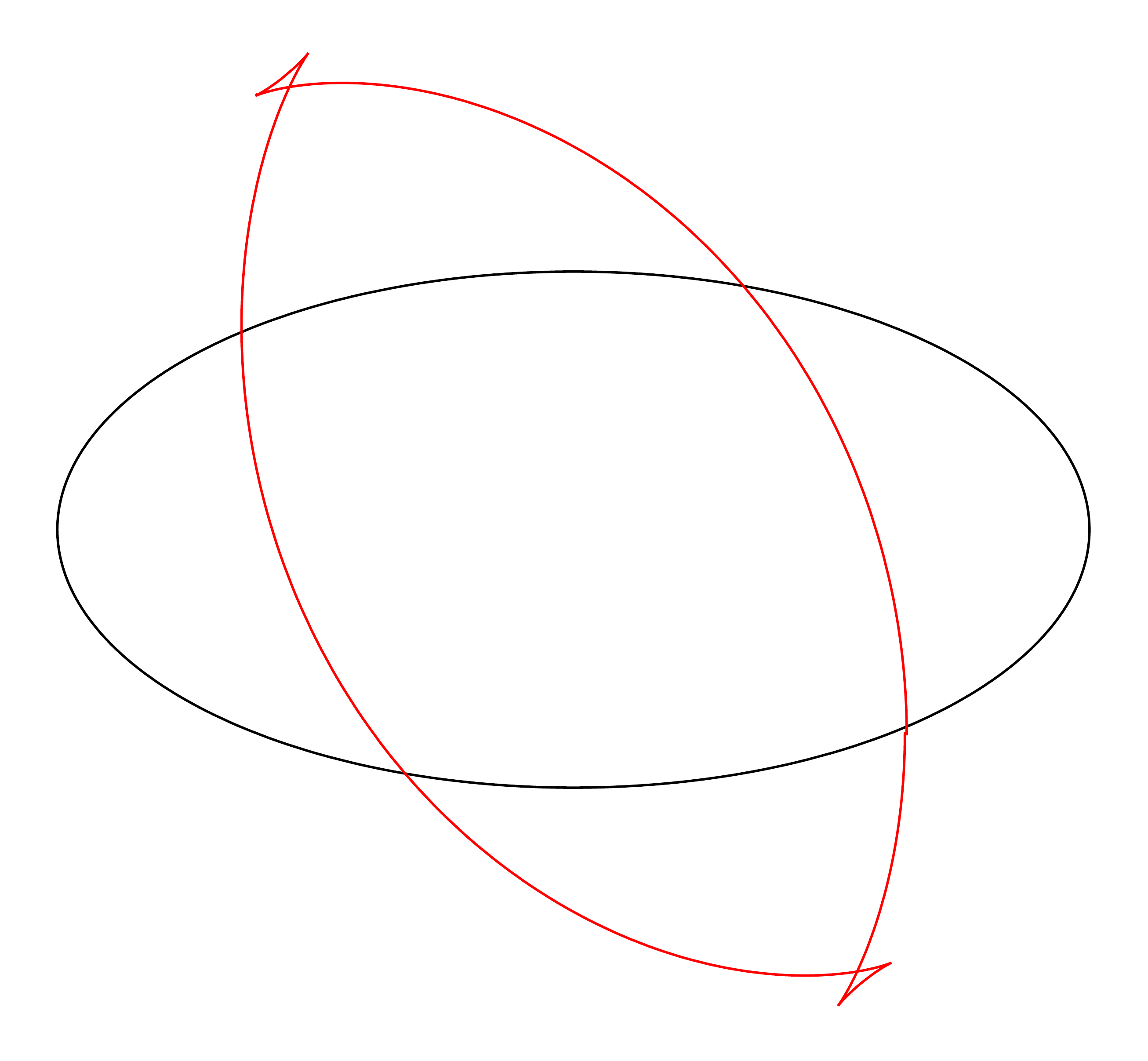}\quad \quad
\includegraphics[width=.4\textwidth]{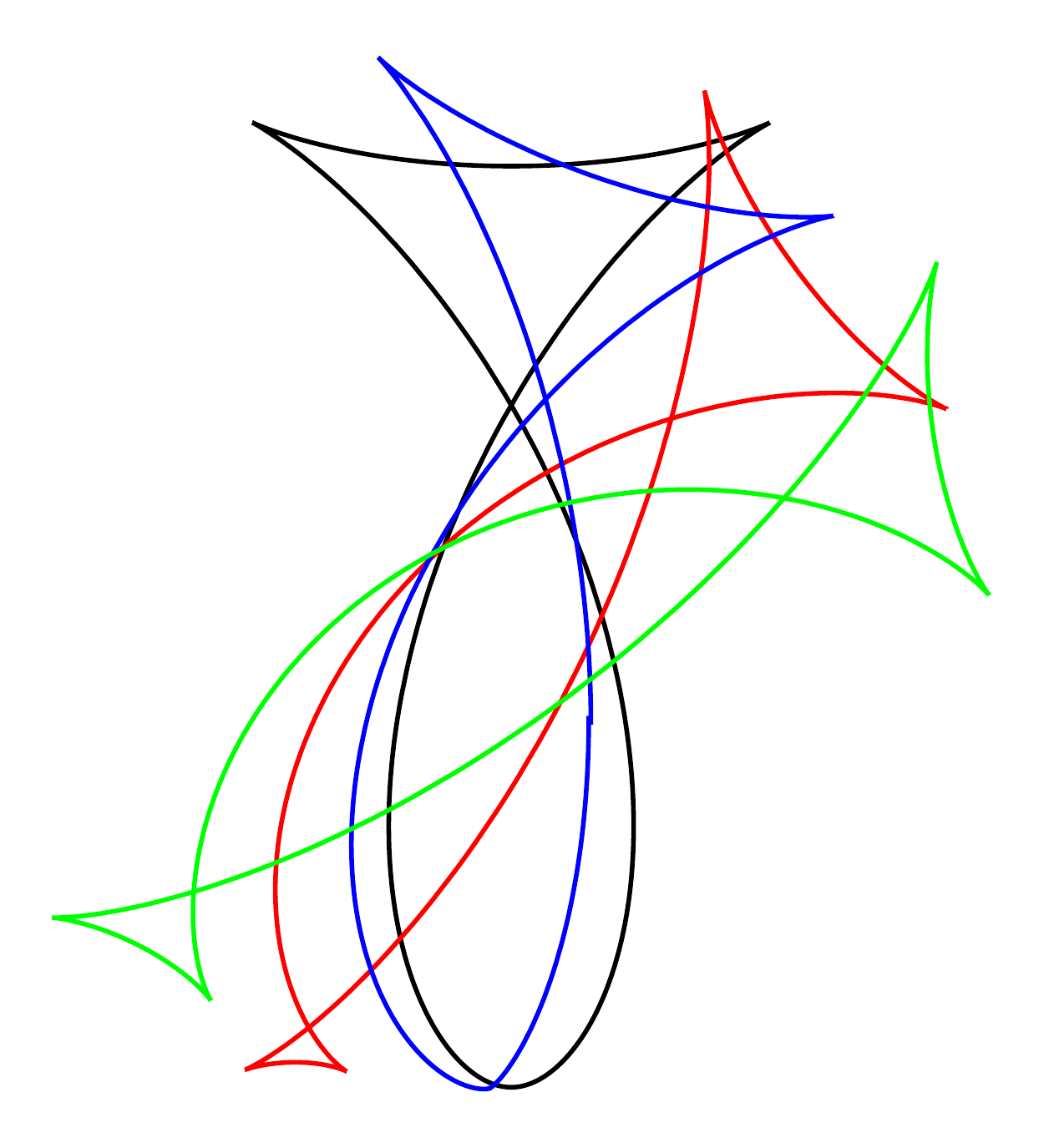}
\caption{Left: the image of an ellipse under the map ${\mathcal M_\alpha}$. Right: the curve with the support function $p(\phi)=e^{2\sin\phi}$ (black) and its images under ${\mathcal M_\alpha}$ for $\alpha=0.5, 0.9$ and $1.2$ (blue, red, and green, respectively).}	
\label{mapellipse}
\end{figure}

Fix an angle $\alpha$ and consider a hedgehog $\G_1$. One constructs its skew-involute $\g$, and then $\G_2$, the skew-evolute of $\g$
with the angle $-\alpha$. We obtain a map  ${\mathcal M_\alpha}={\mathcal E}_{-\alpha}\circ {\mathcal I}_\alpha: \G_1 \mapsto \G_2$. See Figure \ref{mapellipse}.

\begin{figure}[ht]
\centering
\includegraphics[width=.4\textwidth]{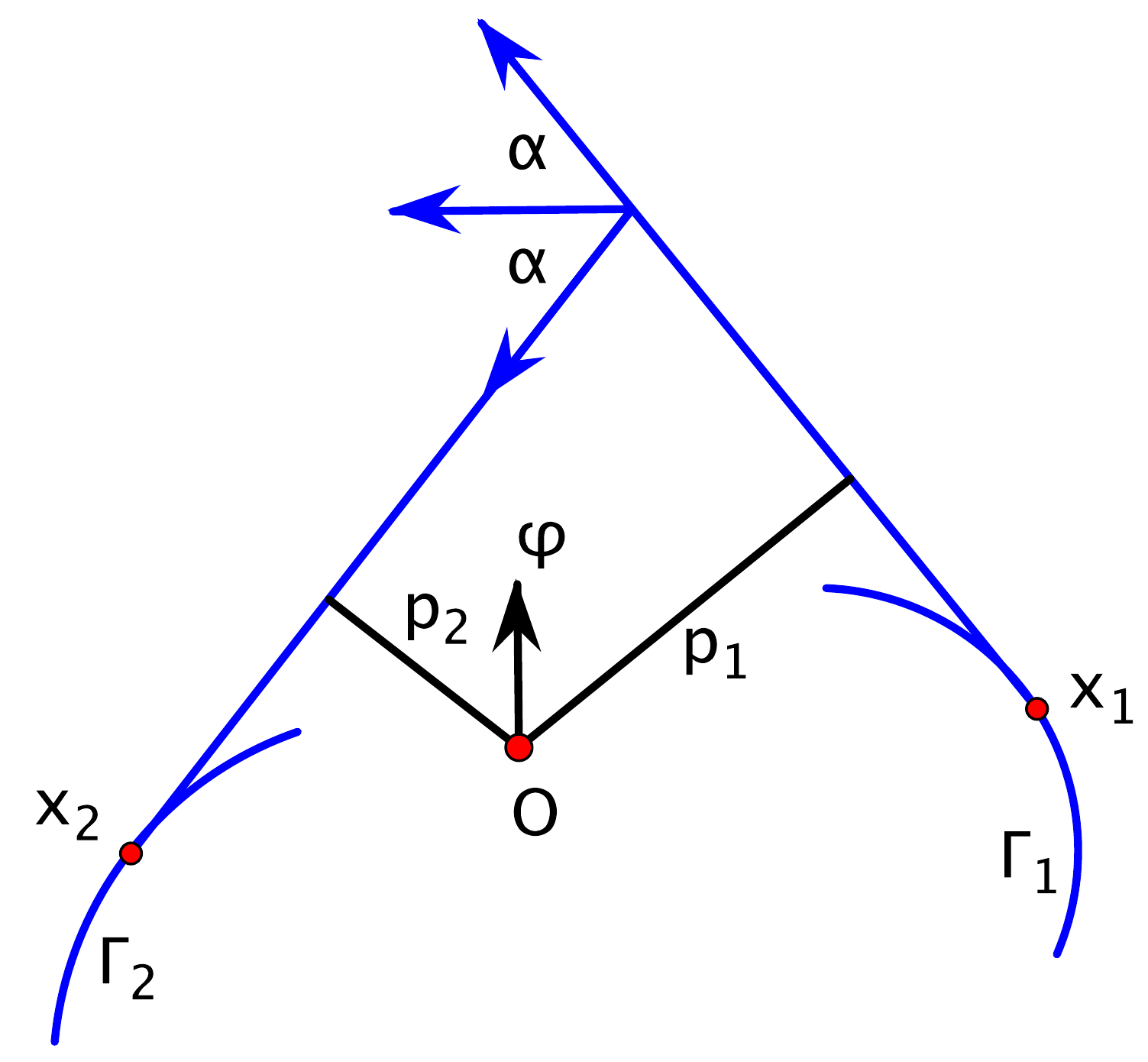}
\caption{The map ${\mathcal M_\alpha}$.}	
\label{map}
\end{figure}

Equivalently, two points, $x_1$ and $x_2$,  traverse the curves $\G_1$ and $\G_2$ in such a way that the  angle between the (co)oriented tangent lines at $x_1$ and $x_2$ is $2\alpha$, and the intersection point of these lines moves in the direction of the bisector between these oriented lines, see Figure \ref{map}.

Let $p_1$ and $p_2$ be the support functions of these curves. The next formula follows from equation (\ref{eq:qtop}):
\begin{equation} \label{eq:map}
p_2'(\phi)-p_2(\phi)\cot\alpha = -p_1'(\phi-2\alpha)-p_1(\phi-2\alpha)\cot\alpha.
\end{equation} 

The next lemma lists some properties of the maps ${\mathcal M_\alpha}$.

\begin{lemma} \label{lm:prop}
1) A curve $\G$ and ${\mathcal M_\alpha}(\G)$ share their Steiner points;\\
2) The maps commute: ${\mathcal M_\alpha} \circ {\mathcal M_\beta} = {\mathcal M_\beta} \circ {\mathcal M_\alpha}$;\\
3) One has ${\mathcal M_\alpha}^{-1} = {\mathcal M_{-\alpha}}$.
\end{lemma}

\proof
The first two properties follow from those of the skew evolute map.
For the third, one has ${\mathcal I_\alpha}={\mathcal E_\alpha}^{-1}$, hence ${\mathcal M_\alpha}={\mathcal E_{-\alpha}} \circ {\mathcal E_\alpha}^{-1}$. It follows that ${\mathcal M_\alpha}^{-1} = {\mathcal E_\alpha} \circ {\mathcal E_{-\alpha}}^{-1} =  {\mathcal M_{-\alpha}}$.
\proofend

The maps ${\mathcal M_\alpha}$ are integrable in the following sense.

\begin{lemma} \label{lm:int}
For every $k$ and every $\alpha$, the sum of the squares of $k$th Fourier coefficients of the support function is preserved by the map ${\mathcal M_\alpha}$:  if 
$$
p(\phi)=\sum_{k\in\Z} a_k e^{ki\phi},\ a_{-k}=\bar a_k,
$$
is a Fourier expansion of the support function of $\G$, then the amplitude $|a_k|$ is an integral of the map ${\mathcal M_\alpha}$ for every $k\ge 0$.
\end{lemma}

\proof
As before, the map preserves the the 2-dimensional spaces of $k$th harmonics. A direct calculation, using equation (\ref{eq:map}), shows that this map is a rotation. More precisely, define the angle $\beta_k$ by $\tan \beta_k=k\tan \alpha$. Then the Fourier coefficients are transformed as follows:
$$
a_k \mapsto a_k e^{2i(\beta_k-\alpha)},
$$
and  the action of ${\mathcal M_\alpha}$ on the space of $k$th harmonics is the rotation by $2(\beta_k-\alpha)$.
\proofend

In particular, hypocycloids evolve by rotations. In this sense, they behave as solitons of the map ${\mathcal M_\alpha}$.

In addition to  the signed length $L(\G)$ and signed area $A(\G)$, let $R(\G)=\int_0^{2\pi} r^2(\phi) d\phi$, where $r(\phi)$ is the curvature radius of the curve $\G$.

\begin{corollary}
One has
$$
L(\G)=L({\mathcal M_\alpha}(\G)),\ A(\G)=A({\mathcal M_\alpha}(\G)),\ {\rm and}\ R(\G)=R({\mathcal M_\alpha}(\G)).
$$
\end{corollary} 

\proof
The first equality directly follows from equation (\ref{eq:map}).

Let
$$
p(\phi)=\sum_{k\in\Z} a_k e^{ki\phi},\ q(\phi)=\sum_{m\in\Z} a_m e^{mi\phi},\ a_{-k}=\bar a_k,\ b_{-m}=\bar b_m,
$$
be the Fourier expansions of two periodic functions. Then
$$
\frac{1}{2\pi} \int_0^{2\pi} p(\phi)q(\phi) d\phi = a_0b_0 + \sum_{k>0} (a_k\bar b_k+\bar a_k b_k).
$$

Let $p(\phi)$ be the support function of $\G$. Then $r(\phi)=p(\phi)+p''(\phi)$, and
$$
A(\G)= \frac{1}{2} \int_0^{2\pi} (p^2(\phi)-p'^2(\phi))d\phi,\quad R(\G)=\int_0^{2\pi} (p(\phi)+p''(\phi))^2 d\phi,
$$
see, e.g, \cite{Sa}. It follows that
$$
\pi A(\G)= a_0^2 +2\sum_{k>0} (1-k^2)|a_k|^2,\ R(\G)= a_0^2 +2\sum_{k>0} (1-k^2)^2 |a_k|^2,
$$
and the result follows from Lemma \ref{lm:int}. 
\proofend

Consider the space of hedgehogs whose support functions are trigonometric polynomials of degree $d$. This space is $2d+1$-dimensional. Fixing the amplitudes of each harmonic, we obtain a space ${\mathcal H}_d$,  a $d$-dimensional torus. If $\G \in {\mathcal H}_d$, then so is ${\mathcal M_\alpha}(\G)$. Geometrically, this space consists of the Minkowski sums of hypocycloids of orders $0,1,\ldots,d$, scaled according to the fixed amplitudes, and each rotated through all angles independently of each other. 

The map ${\mathcal M_\alpha}$ is a rotation of this torus: the $k$th factor $S^1$ is rotated by $2(\beta_k-\alpha)$, where $\beta_k$ are as in the proof of Lemma \ref{lm:int}. For a generic $\alpha$, it is natural to expect the angles $\beta_k$ to be rationally independent.

\begin{conjecture}
For a generic $\alpha$, the orbit of a point is dense in the torus ${\mathcal H}_d$.
\end{conjecture}

\paragraph{Gutkin vs Wegner.}

Circles are  invariant under ${\mathcal M_\alpha}$ for every $\alpha$. {\it Are there other invariant curves?} 

This question is an analog of the following ``bicycle" problem: which curves are in the bicycle correspondence with themselves? This  problem  is equivalent to  Ulam's problem to describe  the bodies that float in equilibrium in all positions (in dimension 2),  problem 19 of The Scottish Book \cite{SC}. 

This problem is not completely solved, but there is a wealth of results in this direction, including constructions of such curves by F. Wegner: these curves are pressurized elastica, and they are solitons of the planar filament equation, a completely integrable partial differential equation of soliton type. See \cite{We1,We2,We3} and \cite{BLPT}.

However, due to linearity, the  problem at hand is considerably simpler, and it was solved by E. Gutkin in the billiard set-up \cite{Gu}. 

Indeed, if ${\mathcal M_\alpha}(\G)=\G$ for a convex curve $\G$, and $\g={\mathcal I_\alpha}(\G)$ is also convex, then $\G$ is a caustic of the billiard inside $\g$, having the special property that the billiard trajectories tangent to $\G$ make angle $\alpha$ with the billiard curve $\g$;  see  also \cite{AA}.

\begin{theorem}[Gutkin] \label{thm:Gu}
A necessary and sufficient condition for such non-circular curves $\G$ to exist is that $k\tan\alpha=\tan(k\alpha)$ for some $k\ge 2$.
\end{theorem} 

\proof To show necessity, set $p_2=p_1=:p$ in (\ref{eq:map}) and rewrite it as
\begin{equation} \label{eq:Gu}
p_2'(\phi+\alpha)\sin\alpha-p_2(\phi+\alpha)\cos\alpha+p_1(\phi-\alpha)\cos\alpha+p_1'(\phi-\alpha)\sin\alpha=0.
\end{equation}
 If
$$
p(\phi)=p_0+\sum_1^\infty a_k \cos (k\phi)+b_k \sin(k\phi),
$$
then equation (\ref{eq:Gu}) implies 
$$
a_k (\sin(k\alpha)\cos\alpha-k\cos(k\alpha)\sin\alpha)=b_k (\sin(k\alpha)\cos\alpha-k\cos(k\alpha)\sin\alpha) =0.
$$
If the curve is not a circle, then $a_k\ne 0$ or $b_k\ne 0$ for some $k\ge 2$, and then
$$
\sin(k\alpha)\cos\alpha=k\cos(k\alpha)\sin\alpha,
$$
as needed.

For sufficiency, one can take a ``fattened" hypocycloid of order $k$, that is, add a sufficiently large constant to the support function of the hypocycloid. This  yields a convex curve having the desired property.
\proofend

\bigskip

{\bf Acknowledgements}. I am grateful to M. Arnold, D. Fuchs, A. Glutsyuk, J. Jer\'onimo-Castro, and I. Izmestiev for  discussions, suggestions, and help. 
I was supported by NSF grant DMS-2005444.


\begin{thebibliography}{99}

\bibitem{AA} V. Aguilar-Arteaga, R. Ayala-Figueroa, I. Gonz\'alez-Garc\'ia, J. Jer\'onimo-Castro. {\it On evolutoids of planar convex curves II. } Aequationes Math. {\bf 89} (2015), 1433--1447.  

\bibitem{AM} T. Apostol, M. Mnatsakanian. {\it Tanvolutes: generalized involutes.}  Amer. Math. Monthly {\bf 117} (2010), 701--713. 

\bibitem{AFITT} M. Arnold, D. Fuchs, I. Izmestiev, S. Tabachnikov, E. Tsukerman.
{\it Iterating evolutes and involutes}. Discr. Comp. Geom. {\bf 58} (2017), 80--143.

\bibitem{BLPT} G. Bor, M. Levi, R. Perline, S. Tabachnikov. {\it Tire tracks and integrable curve evolution.} Int. Math. Res. Notes, 2020, no. 9, 2698--2768.

\bibitem{FLT} R. Foote, M. Levi, S. Tabachnikov. {\it Tractrices, bicycle tire tracks, hatchet planimeters, and a 100-year-old conjecture.} Amer. Math. Monthly {\bf 120} (2013), 199--216.

\bibitem{GW} P. Giblin, J. Warder. {\it Evolving evolutoids.} Amer. Math. Monthly {\bf 121} (2014), 871--889. 

\bibitem{Gu} E. Gutkin. {\it Capillary floating and the billiard ball problem.} J. Math. Fluid Mech. {\bf 14} (2012), 362--382.

\bibitem{Ha} M. Hamman. {\it  A note on ovals and their evolutoids}. Beitr.  Algebra Geom. {\bf 50} (2009), 433--441.

\bibitem{HPZ} S. Howe,  M. Pancia, V. Zakharevich. {\it Isoperimetric inequalities for wave fronts and a generalization of Menzin's conjecture for bicycle monodromy on surfaces of constant curvature.} Adv. Geom. {\bf 11} (2011), 273--292.

\bibitem{JC} J. Jer\'onimo-Castro. {\it On evolutoids of planar convex curves.}  Aequationes Math. {\bf 88} (2014), 97--103. 

\bibitem{LT} M. Levi, S. Tabachnikov {\it On bicycle tire tracks geometry, hatchet planimeter, Menzin?s conjecture, and oscillation of unicycle tracks}.
Experiment. Math. {\bf 18} (2009), 173--186.

\bibitem{PW} G. P\'olya, N. Wiener. {\it On the oscillation of the derivatives of a periodic function.} Trans. Amer. Math. Soc. {\bf 52} (1942), 249--256.

\bibitem{Re} R.-A. F. de R\'eaumur. {\it  Methode g\'en\'erale pour d\'eterminer le point d'intersection de deux lignes droites
infiment proches, qui rencontrent une courbe quelconques vers le m\^eme c\^ot\'e sous des angles \'egaux moindres, ou plus grandes qu'un droit, et pour conno\^itre la nature de la courbe d\'ecrite par une infinit\'e de tels points d'intersection.} Histoire de l'Acad\'emie Royale des Sciences (1709).

\bibitem{RGS} D. Reznik, R. Garcia, H. Stachel. {\it Area-invariant pedal-like curves derived from the ellipse.} Beitr. Algebra Geom. {\bf 63} (2022), 359--377.

\bibitem{Sa} L. Santal\'o. Integral geometry and geometric probability. Second edition. Cambridge Univ. Press, Cambridge, 2004. 

\bibitem{SC} The Scottish Book. Mathematics from the Scottish Caf\'e with selected problems from the new Scottish Book. Second edition.  Edited by R. Daniel Mauldin. Birkh\"auser/Springer, Cham, 2015. 



\bibitem{We1} F. Wegner. {\it Floating Bodies of Equilibrium in 2D, the Tire Track Problem and Electrons in a Parabolic Magnetic Field}. arXiv:physics/0701241. 

\bibitem{We2} F. Wegner. {\it From Elastica to Floating Bodies of Equilibrium}. arXiv:1909.12596. 

\bibitem{We3} F. Wegner. {\it Three and more Problems -- One Solution}. \url{https://www.thphys.uni-heidelberg.de/~wegner/Fl2mvs/Movies.html}

\bibitem{Wu} W. Wunderlich. {\it \"Uber die Evolutoiden der Ellipse.} Elem. Math. {\bf 10} (1955), 37--40.

\bibitem{Z} M. Zwierzy\'nski. {\it The singular evolutoids set and the extended evolutoids front.} Aequationes Math. {\bf 96} (2022), 849--866.

\end{thebibliography}
\end{document}